\documentclass[10pt,conference]{IEEEtran} 

\IEEEoverridecommandlockouts    

\title{\LARGE \bf
Discounted MPC and infinite-horizon optimal control under plant-model mismatch: Stability and suboptimality 
}
\author{Robert H. Moldenhauer, Karl Worthmann, Romain Postoyan, Dragan Ne\v{s}i\'{c} and Mathieu Granzotto % <-this % stops a space
\thanks{Work supported by the Australian Research Council under the Discovery Grant DP250100300, the ANR grant OLYMPIA ANR-23-CE48-0006, and the DFG grant 535860958 within the research unit ALeSCo.}% <-this % stops a space
\thanks{R. H. Moldenhauer, D. Ne\v{s}i\'{c} and M. Granzotto are with the Department of Electrical and Electronic Engineering, University of Melbourne, Parkville, VIC 3010, Australia (e-mail: moldenhauer.r@student.unimelb.edu.au, dnesic@unimelb.edu.au, mathieu.granzotto@unimelb.edu.au).}
\thanks{R. H. Moldenhauer and R. Postoyan are with the Université de Lorraine, CNRS, CRAN, F-54000 Nancy, France (e-mail: \{name.surname\}@univ-lorraine.fr).}
\thanks{K. Worthmann is with the Optimization-based Control Group, %Institute of Mathematics, 
TU Ilmenau, 98693 Ilmenau, Germany (e-mail: karl.worthmann@tu-ilmenau.de).}
}

\usepackage{cite}
\usepackage{amsmath,amssymb,amsfonts}
\usepackage{algorithmic}
\usepackage{graphicx}
\usepackage{algorithm,algorithmic}
\usepackage{hyperref}
\usepackage{textcomp}
\usepackage{url}
\usepackage{soul}
\usepackage{xcolor}
\usepackage{mathtools}

\usepackage{svg,booktabs}
\let\labelindent\relax % Added by Karl to resolve a conflict
\usepackage{enumitem}
\usepackage{theoremref}
\usepackage{amssymb}

\usepackage{amsthm}
\usepackage{hyperref}
\usepackage{nicefrac}
\usepackage{mathtools}
\hypersetup{hidelinks}
\theoremstyle{thm}
\newtheorem{thm}{Theorem}
\newtheorem{prop}{Proposition}
\newtheorem{corollary}{Corollary}
\newtheorem{remark}{Remark} 
\newtheorem{lem}{Lemma}

\newtheorem{standingassm}{Standing Assumption}

\newcommand\numberthis{\addtocounter{equation}{1}\tag{\theequation}}
\newcommand{\R}{\mathbb{R}}
\newcommand{\N}{\mathbb{N}}

\newcommand{\mc}{\mathcal}

\newcommand{\okappa}{\overline\kappa}

\newcommand{\lmax}{{\overline\lambda}}
\newcommand{\lmin}{{\underline\lambda}}

\usepackage{tikz}
\usetikzlibrary{matrix, arrows.meta}

\providecommand{\phantomsection}{}
\makeatletter
\newcommand{\taglabel}[2]{%
    \text{(#1)}%
    \phantomsection
    \protected@edef\@currentlabel{#1}%
    \label{#2}%
}
\makeatother

\newlength{\rdpw}
\newlength{\rdptmp}
\newcommand{\rdpmax}[1]{%
    \settowidth{\rdptmp}{#1}%
    \ifdim\rdptmp>\rdpw \setlength{\rdpw}{\rdptmp}\fi
}

\newcommand{\gnm}{{\gamma,N-1}}
\newcommand{\gn}{{\gamma,N}}
\newcommand{\gi}{{\gamma,\infty}}
\newcommand{\gnn}{{\gamma,N_0}}

\newcommand{\Ni}{\N\cup\{\infty\}}

\newcommand{\oi}{{1,\infty}}

\newcommand{\oll}{{\overline\lambda}}

\newcommand{\No}{{N_0}}

\DeclareMathOperator*{\argmin}{arg\,min}

\newcommand{\agn}{{\alpha_\gn}}
\newcommand{\aggn}{{\alpha_{\underline{\smash{\gamma}},N}}}
\newcommand{\agnn}{{\alpha_{\gamma,\underline N}}}
\newcommand{\aggnn}{{\alpha_{\underline{\smash{\gamma}},\underline N}}}

\usepackage{ushort}

\begin{document}

\maketitle

% REQUIRED
\begin{abstract}
    We study the robustness to plant-model mismatch of nonlinear systems in closed loop with (discounted) optimal controllers designed for a nominal model, with quadratic stage costs across finite and infinite horizons. By combining Lyapunov analysis with a unifying relaxed dynamic programming framework, we establish %explicit 
    conditions that guarantee robust stability alongside improved suboptimality bounds compared to existing literature. A key step in our analysis is proving that stability counteracts time-accumulated error growth due to plant-model mismatch, yielding robustness bounds that hold uniformly across all prediction horizons.
\end{abstract}

% REQUIRED
%\begin{keywords}
%Predictive control for nonlinear systems, stability of nonlinear systems, optimal control
%\end{keywords}

\section{Introduction}
    In optimal and model predictive control (MPC), the study of robustness has %The study of robustness Robustness of optimal and model predictive control (MPC) is a fundamental question with 
    a rich history: R. E. Kalman showed in 1964 that the linear quadratic regulator (LQR) possesses guaranteed phase and gain margins \cite{kalman_when_1964}.
    These margins stem from dissipativity and were generalized to control-affine nonlinear systems, see, e.g., \cite{sepulchre_constructive_1997}. %and references therein.
    However, such guarantees rely on specifically structured uncertainties, %specific uncertainty structures, 
    beyond which optimal controllers can become %exhibit 
    fragile w.r.t.\ robustness \cite{doyle_guaranteed_1978,sussmann_peaking_1991, grimm_examples_2004}.
    While frameworks such as robust adaptive dynamic programming have been developed \cite{jiang2017robust}, many open questions remain regarding robustness under plant-model mismatch, especially %those arising 
    in data-driven and learning-based %optimal 
    control.
    A particularly successful approach to robustness analysis relies on regular Lyapunov functions, %with regularity properties (e.g., continuity, local Lipschitz). 
    %If such a Lyapunov function exists for the unperturbed model, it 
    which may also serve as a candidate Lyapunov function for the perturbed system. %to establish stability under broad classes of perturbations.
    This approach %framework 
    has been applied %express plant-model mismatch 
    across numerous %practical 
    scenarios, including parametric uncertainty \cite{khalil2002nonlinear}, sampling and discretization~\cite{nesic_sufficient_1999}, measurement errors and state estimation \cite{messina_discrete-time_2005}, %unmodeled dynamics in simplified design models \cite{ioannou_instability_1984}, and models identified from data
    and data-driven modeling \cite{dean_sample_2020, berberich_data-driven_2021, BoldGrun25, SchiWort25}.

    In receding horizon %optimal 
    control, the optimal value function often serves as a Lyapunov function candidate for the nominal %unperturbed 
    closed-loop system using %. Specifically, %the dynamic programming (DP) equation or a 
    relaxed dynamic programming (DP), see %inequality 
    %can be used to establish %establishes that 
    %the value function satisfies Lyapunov conditions 
    \cite{GrunRant08,grune_analysis_2009}. %,grune_economic_2013} under mild assumptions. 
    This approach extends to discounted costs \cite{postoyan_stability_2017, granzotto_finite-horizon_2021, moldenhauer2026stability} and can provide horizon-dependent near-optimal performance guarantees.
    When the value function satisfies certain continuity properties, this framework yields explicit robustness guarantees for MPC \cite{grimm2007nominally, BoldGrun25, SchiWort25, liu_certainty-equivalence_2026, kuntz_beyond_2026}.
    However, existing robustness bounds \cite{SchiWort25,liu_certainty-equivalence_2026} generally deteriorate as the prediction horizon grows and fail entirely for infinite horizons. Moreover, these bounds are often conservative, which limits their practical applicability, and we are not aware of any result applicable to discounted costs.

To overcome these limitations, we establish robustness %robust stability and near-optimality
guarantees for nonlinear optimal control over finite horizons (MPC) and infinite horizons with discounted and undiscounted quadratic stage costs. Controllers are designed on a \emph{surrogate model} and deployed on the true plant. 
We develop a unifying %single 
relaxed DP framework encompassing %that treats the aspects of 
plant-model mismatch, horizon length, and discounting. %jointly.
It yields explicit bounds on these three ingredients for exponential stability, near-optimality, and continuous dependence of the value function on the model.
%These bounds reveal that all three aspects share a common robustness margin, which discounting enlarges by smoothing the value function. 
In contrast to the horizon-dependent results in the literature, our %mismatch 
bounds hold uniformly in the horizon length. %over all finite and infinite horizons.
This is important when a long horizon is desirable in MPC. 
To this end, we exploit %Our results reveal that 
stability to counteract %counteracts 
the accumulation of disturbances along the optimization horizon. %to deterioration a long-horizon cost, which therefore does not harm robustness.
%and our results show that this does harm robustness, because stability counteracts the accumulation of disturbances.
We show on an example that the derived bounds on the plant-model mismatch improve upon those in \cite{liu_certainty-equivalence_2026} even for short horizons.
Furthermore, we demonstrate that discounting the stage cost can also increase these bounds.
\textbf{Notation}. $\R$ ($\R_{\geq0}$) denotes the set of (nonnegative) real numbers and $\N$ ($\N_{0}$) the set of positive (nonnegative) integers. 
%The empty set is denoted  $\varnothing$. 
Given a real symmetric, positive definite matrix $P$, we write $\|x\|_P := \sqrt{x^\top Px}$ for any $x\in\R^n$ and denote its largest and smallest eigenvalue by $\lmax_P$ and $\lmin_P$, respectively. %Given $a\in\R^n$, $n\in\N$, $\operatorname{diag}(a)$ is the $n\times n$ diagonal matrix, whose diagonal elements form the vector $a$. 
A function $\alpha:[0,\infty)\to[0,\infty)$ is of class-$\mathcal{K}$, denoted by $\alpha\in\mathcal{K}$, if it is continuous, zero at zero, and strictly increasing.

\section{Problem statement}\label{sect:problem-statement}
First, we introduce the plant and the model dynamics in Section~\ref{s:2A}. 
Second, we formulate the optimal control problem (OCP) and introduce the value function in Section~\ref{subsect:optimal-control-problem}, before presenting the problem statement in Section~\ref{s:2C}.

\subsection{Plant, surrogate model, and plant-model mismatch}\label{s:2A}
%\subsection{Plant and surrogate models}\label{s:2A}
Consider the discrete-time system
\begin{align}\label{eq:system-dynamics}
    x^+ = g(x,u)
\end{align}
with state $x\in\R^n$, control $u\in\mathbb U\subseteq\R^m$, and successor state $x^+ \in \mathbb{R}^n$, where %$n,m\in\N$. 
%Let 
the set $\mathbb U$ %of feasible control 
satisfies the following.
\begin{standingassm}[\textbf{SA1}]\label{assm:1:setU}
    The set~$\mathbb U$ of admissible controls is closed and contains the origin.
    %in its interior.
\end{standingassm}
Let $g:\R^n\times\mathbb U\to\R^n$ satisfy $g(0,0)=0$, i.e., the origin is an equilibrium.
In practice, the function $g$ may be unknown or too complex for control design.
Hence, we consider the scenario where control inputs to \eqref{eq:system-dynamics} are designed using a \emph{surrogate model}
\begin{align}
    x^+ = f(x,u) \label{eq:surrogate}
\end{align}
with a known function $f:\R^n\times\mathbb U\to\R^n$.
This function is typically obtained by modeling, system identification~\cite{schoukens2019nonlinear}, or data-driven techniques~\cite{kutz2016dynamic}.
Denote by $\varphi^g_x(k,\mathbf u_k)$ and $\varphi^f_x(k,\mathbf u_k)$ the states at time $k\in\N_0$ of plant~\eqref{eq:system-dynamics} and model~\eqref{eq:surrogate}, respectively, where $x\in\R^n$ and $\mathbf u_k = (u_0,\dots,u_{k-1})\in\mathbb U^k$ denote the initial state at time~$0$ and the control sequence, respectively. 
We assume that $f$ is continuous in $u$ and Lipschitz in $x$ uniformly in $u$.
%, as stated below.
\begin{standingassm}[\textbf{SA2}]\label{assm:Lipschitz}
    The function $f$ is continuous and there exists a constant $L\geq0$ such that, for all $x_1,x_2 \in \R^n$ and $u\in\mathbb U$,
\begin{align}\label{eq:L-Lipschitz}%\tag{L-Lipschitz}
    |f(x_1,u) - f(x_2,u)| \leq L |x_1 - x_2|.
\end{align}
%holds .
\end{standingassm}

We consider the \emph{proportional (plant-model) mismatch} on a given set $\mc S\subseteq\R^n$, defined as
\begin{align}
    |f-g|_{\mc S} := \sup_{(x,u)\in(\mc S\times\mathbb U)\setminus\{(0,0)\}}%{x\in\mc S\setminus\{0\}, u\in\mathbb U}
    \frac{|f(x,u)-g(x,u)|}{|x|+|u|},
    %\inf\big\{\overline p&\geq0:|f(x,u)-g(x,u)|\nonumber\\&~~~~\leq \overline p(|x|+|u|)~\forall\,x\in\mc S, u\in\mathbb U\big\}
    \label{eq:defn:mismatch}
\end{align} 
with the following standing assumption.
\begin{standingassm}[\textbf{SA3}]\label{assm:3:mismatch}
    The set $\mc S\subseteq\R^n$ contains the origin~$0$, and it holds that $f(0,0)=0$ and $|f-g|_{\mc S}<\infty$.
\end{standingassm}
The proportional mismatch \eqref{eq:defn:mismatch} is a form of \emph{parametric uncertainty}, which means that the allowed magnitude of uncertainty $|f(x,u)-g(x,u)|$ depends on the state $x$ (and control $u$).
In particular, we require that the model is exact at $0$ (i.e., $f(0,0)=g(0,0)=0$), and uncertainty may increase linearly with $|x|+|u|$.
Such a proportional mismatch was considered in, e.g., \cite{BoldGrun25,liu_certainty-equivalence_2026,SchiWort25,kuntz_beyond_2026}.
The set~$\mc S$ is to be interpreted as a region in the state space where a bound on the mismatch is known (and sufficiently small).
For example, $\mc S$ can be the set where certain modeling assumptions are reasonable, e.g., a spring behaving approximately linearly.
Note that, even when $g$ is not exactly known, bounds for $|f-g|_{\mc S}$ may still be found, see, e.g., \cite{grune2003optimization} (discretization errors due to numerical approximations) or \cite{PaulKoch21,PaulGram24} (data-driven surrogates using neural networks).
For data-driven identification techniques, $\mc S$ represents a region in which sufficient data are available.
It is shown in~\cite{BoldPhil25,SchiWort25} that, if the plant~$g$ is Lipschitz in $x$ uniformly in $u$ and control affine, then kernel extended dynamic mode decomposition (kEDMD) can be used to obtain data-driven surrogates~$f$ with arbitrarily small proportional mismatch $|f-g|_{\mc S}$ on any compact set~$\mc S$ on which sufficient data are available.

\begin{remark}\label{rem:local_assumptions}
    We assume \eqref{eq:L-Lipschitz}  globally, i.e., for all $x_1,x_2\in\R^n$, for simplicity.
    If it is satisfied only locally, e.g., on $\mc S$, then $f$ may be redefined outside of $\mc S$ to satisfy \eqref{eq:L-Lipschitz} globally, without changing $|f-g|_{\mc S}$.
\end{remark}

\subsection{Optimal control problem}\label{subsect:optimal-control-problem}

We are interested in the setting where %optimal 
controls are computed using %along solutions to 
the surrogate model~\eqref{eq:surrogate}, but are applied to the plant~\eqref{eq:system-dynamics}.
Therefore, for $N \in \mathbb{N} \cup \{\infty\}$, %first 
consider the optimal control problem for any given state $x\in\R^n$
\begin{align}
    \inf_{\mathbf{u}_N\in\mathbb U^N} J_\gn^f(x,\mathbf u_N)
    \label{eq:OCP}
\end{align}
with the cost function $J_\gn^f:\R^n\times\mathbb U^N\to\R_{\geq0}\cup\{\infty\}$ defined for $x\in\R^n$ and $\mathbf u_N=(u_0,\dots,u_{N-1})\in\mathbb U^N$ as
\begin{align}\label{eq:cost-J}
    J_\gn^f(x,\mathbf u_N) := \sum_{k=0}^{N-1}\gamma^k\ell(\varphi^f_x(k,\mathbf{u}_k),u_k).
\end{align}
The stage cost $\ell:\R^n\times\mathbb U\to\R_{\geq0}$ is quadratic and given by
\begin{align}
    \ell(x,u) := x^\top Qx + u^\top Ru \label{eq:stage_cost}
\end{align}
with weighting %fixed 
matrices $Q\in\R^{n\times n}$ and $R\in\R^{m\times m}$.
We assume the following condition. 
\begin{standingassm}[\textbf{SA4}]\label{assm:4:QR}
    The matrices $Q \in \R^{n\times n}$ and $R \in \R^{m\times m}$ in \eqref{eq:stage_cost} are symmetric and positive definite.
\end{standingassm}
The number~$N \in \Ni$ in~\eqref{eq:cost-J}, also termed prediction or optimization horizon, %shall take values in $\Ni$ thereby allowing 
enables us to consider both finite- and infinite-horizon cost functions in a unified way. The constant  $\gamma\in(0,1]$ is the discount factor.
For $\gamma < 1$, future costs are weighted less with an exponential discount.
This is commonly used in dynamic programming as it renders the Bellman operator a contraction under mild conditions \cite{bertsekas_dynamic_2012}.
Further, we define
the \emph{(optimal) value function} $V_\gn^f:\R^n\to\R\cup\{\infty\}$ as %, for any $x\in\R^n$,
\begin{align}
    V_\gn^f(x) := \inf_{\mathbf u_N\in\mathbb U^N}J_\gn^f(x,\mathbf u_N),
    \label{eq:defn:value_function}
\end{align}
and $V_{\gamma,0}^f:\equiv0$.
Similarly, the value function for \eqref{eq:system-dynamics} is denoted by $V_\gn^g(x) := \inf_{\mathbf u_N\in\mathbb U^N} \sum_{k=0}^{N-1}\gamma^k\ell(\varphi^g_x(k,\mathbf{u}_k),u_k)$. % denotes the value function for the plant \eqref{eq:system-dynamics}. %we define $J_\gn^g$ and $V_\gn^g$ for the plant~\eqref{eq:system-dynamics}.
We require the following controllability property of the surrogate 
model~\eqref{eq:surrogate} in terms of its %quadratic 
optimal cost.
\begin{standingassm}[\textbf{SA5}]\label{assm:cost_controllable}
    There exists $B\geq1$ such that, for all $x\in\R^n$,
    \begin{align}\label{eq:B-cost-controllable}%\tag{B-cost-controllable}
        V_{1,\infty}^f(x) \leq B\|x\|_Q^2.
    \end{align} 
\end{standingassm}
It was shown in \cite[Lemma 3.2]{grune_analysis_2009}, see also \cite[Lemma 6.8]{grune_nonlinear_2017} and \cite[Lemma 1]{postoyan_stability_2017}, that SA\ref{assm:cost_controllable} holds if the stage cost is uniformly globally exponentially controllable to zero\footnote{In the sense that there exist $M>0$ and decay rate $\lambda>0$ %, where $\lambda$ is the \emph{decrease rate}, 
such that, for any $x \in \R^n$, there exists an infinite-length control sequence $\mathbf u_\infty\in\mathbb U^\infty$ satisfying $\ell(\varphi^f_x(k,\mathbf u_k), u_k) \leq M|x|^2e^{-\lambda k}$ for all $k\in\N_0$.} along solutions to the surrogate model \eqref{eq:surrogate}.
The converse implication also holds, as it is a special case of our main result that, under SA\ref{assm:cost_controllable}, optimal controls for $\gamma=1$, $N=\infty$ and $f=g$ achieve exponential stability with exponentially decaying controls. 
Hence, SA\ref{assm:cost_controllable} is a mild assumption when the objective is to ensure %achieve 
exponential stability; see \cite[A2]{tuna2006shorter} and~\cite{GrunRant08} for similar assumptions.
Note that Remark~\ref{rem:local_assumptions} similarly applies to SA\ref{assm:cost_controllable}.

\begin{remark}[\textbf{Comparison to \cite{liu_certainty-equivalence_2026}}]\label{rem:comparison_liu}
    Similar assumptions were used in \cite{liu_certainty-equivalence_2026} to guarantee robustness and suboptimality under plant-model mismatch. %\cite{SchiWort25} (which considered question 1 of Section \ref{s:2C}) and  \cite{liu_certainty-equivalence_2026} (which considered all questions of Section \ref{s:2C}).
    However, the results in \cite{liu_certainty-equivalence_2026} require additional regularity assumptions
    for the optimal inputs (see \cite[Assumption 5]{liu_certainty-equivalence_2026}). 
    Our approach relies solely on regularity of the value function and does not need such assumptions.
    %We do without such assumption, relying solely on the regularity of the value function. %Our methodology requires no such assumptions.
\end{remark}

\begin{remark}[\textbf{Relaxing SA\ref{assm:4:QR}}]
For clarity, we restrict ourselves to quadratic stage cost with $Q,R\succ0$, but the results can be adapted to non-quadratic stage costs with quadratic bounds \cite{liu_certainty-equivalence_2026}, to positive semidefinite $R$ if the mismatch in \eqref{eq:defn:mismatch} is measured relative to $|x|$ instead of $|x|+|u|$, and to positive semidefinite $Q$ under a suitable detectability condition \cite{postoyan_stability_2017, granzotto_finite-horizon_2021}.
\end{remark}

\begin{remark}[\textbf{Further relaxations}]
    The assumptions can further be relaxed to guarantee notions weaker than exponential stability.
    Allowing mismatch at the origin yields practical exponential stability, while merely continuous dynamics and milder controllability yield asymptotic stability under certain state-dependent mismatch, and semi-global practical stability under persistent mismatch \cite{grimm2007nominally, nesic_sufficient_1999, BoldGrun25, postoyan_stability_2017, granzotto_finite-horizon_2021, moldenhauer2026stability}.
\end{remark}

%Under SA1 and SA2, if finite costs are achievable and $f$ is continuous, then 
Under our assumptions, the OCP \eqref{eq:OCP} always has a solution and the infimum in \eqref{eq:defn:value_function} can be replaced by a minimum. %as a minimum. 
This, along with the Bellman equation \cite{bertsekas_dynamic_2012}, is stated in the next  proposition, whose proof directly follows from %by application of 
\cite[Theorem~2]{keerthi2003existence} and is therefore omitted (see also \cite{bertsekasValuePolicyIterations2017}).
\begin{prop}[\textbf{%Existence of o
Optimal controls \& Bellman equation}]\label{prop:1.bellman}
    %Let system \eqref{eq:surrogate} be given and consider the OCP~\eqref{eq:OCP}
    Consider the system \eqref{eq:surrogate} and the OCP~\eqref{eq:OCP} 
    with discount factor $\gamma\in(0,1]$ and horizon %length 
    $N\in\N\cup\{\infty\}$.
    Then, for every state $x\in\R^n$, %satisfying $V_\gn^f(x)<\infty$,
    there exists a control sequence $\mathbf u_N^\star\in\mathbb U^N$ such that
    % %\begin{subequations}\label{eq:value-function}
    % %\begin{align}
    %     $V_\gn^f(x) = J_\gn^f(x,\mathbf u_N^\star)%\label{eq:defn:value_function}\\
    %     %= \min_{\mathbf u_N\in\mathbb U^N}J_\gn^f(x,\mathbf u_N)\\
    %     = \min_{u\in\mathbb U}\big(\ell(x,u) + \gamma  V_\gnm^f(f(x,u))\big).$%\label{eq:prop:Bellman}
    % %\end{align}
    % %\end{subequations}
    \begin{align}
        V_\gn^f(x) &= J_\gn^f(x,\mathbf u_N^\star)\nonumber\\
        &= \min_{u\in\mathbb U} \left(\ell(x,u) + \gamma  V_\gnm^f(f(x,u))\right).
    \end{align}
\end{prop}
We define the \emph{set-valued optimal %feedback 
policy} $\mc U_\gn^f:\R^n\rightrightarrows\mathbb U$,
\begin{align}
    \mc U_\gn^f(x):=\argmin_{u\in\mathbb{U}}\left\{\ell(x,u) + \gamma V_\gnm^f(f(x,u))\right\}\label{eq:UgammaN}
\end{align}
%if $V_\gn^f(x)<\infty$, and $\mc U_\gn^f(x):=\varnothing$ otherwise.
for all $x\in\R^n$.
The set $\mc U_\gn^f(x)$ contains the first element of each optimal sequence for the OCP~\eqref{eq:OCP} at state $x$. %. provided $V_\gn^f(x)<\infty$. 
In particular, the Bellman equation 
\begin{align}
\gamma V_\gnm^f(f(x,u)) - V_\gn^f(x) = -\ell(x,u) \tag{DP}\label{eq:DP}    
\end{align}
holds for all $x\in\R^n$ and $u\in\mc U_\gn^f(x)$.

\subsection{Problem formulation}
\label{s:2C}

%So far, we considered the OCP \eqref{eq:OCP} only on the surrogate model \eqref{eq:surrogate}.
In practice, the optimal inputs in \eqref{eq:UgammaN} %the inputs solving the OCP~\eqref{eq:OCP} %optimal inputs obtained by solving \eqref{eq:OCP} 
are not applied to the surrogate model~\eqref{eq:surrogate}, but to the real plant~\eqref{eq:system-dynamics}.
%Hence, we aim to study the difference inclusion
This leads to the closed-loop system given by the difference inclusion
\begin{align}
    x^+ \in g\big(x,\mc U_\gn^f(x)\big) := \big\{g(x,u)~|~u\in\mc U_\gn^f(x)\big\}.\label{eq:closed_loop}
\end{align}
%This is the closed loop that arises when a (discounted) optimal MPC controller is designed using the surrogate model~\eqref{eq:surrogate} in the optimization step, but is applied to the actual plant~\eqref{eq:system-dynamics}.
We are interested in the following three questions.

\textit{1) Exponential stability %of the origin w.r.t.\ the 
of the closed loop~\eqref{eq:closed_loop}}, i.e., 
%The first is about the \emph{exponential stability} of the closed loop \eqref{eq:closed_loop}. That is, 
whether computing optimal inputs minimizing~\eqref{eq:cost-J} w.r.t.\ %the objective function~\eqref{eq:cost-J} along solutions of 
the surrogate model~\eqref{eq:surrogate} leads to exponentially stabilizing policies for the %actual 
plant~\eqref{eq:system-dynamics} at the same equilibrium (the origin).

\textit{2) Suboptimality}, i.e., quantifying %how good 
the performance %remains 
under plant-model mismatch compared to the optimal costs computed for $f$ and $g$. 
To this end, we define the %, for $\gamma\in(0,1]$ and $N\in\Ni$, 
\emph{performance indices} evaluating the (worst-case) discounted and undiscounted cost along all solutions to \eqref{eq:closed_loop}, i.e.,
\begin{align}
    &\overline J_\gi^g (x,\mc U_{\gamma,N}^f)\label{eq:performance1}\\\nonumber
    &:=\sup\left\{J_\gi^g(x,\mathbf u_\infty) \,\big|\, u_{k}\in\mc U_{\gamma,N}^f(\varphi^g_x(k,\mathbf u_k))\,\forall k\in\N_0\right\}\\
    &\overline J_\oi^g (x,\mc U_{\gamma,N}^f)\label{eq:performance2}\\\nonumber
    &:=\sup\left\{J_\oi^g(x,\mathbf u_\infty) \,\big|\, u_{k}\in\mc U_{\gamma,N}^f(\varphi^g_x(k,\mathbf u_k))\,\forall k\in\N_0\right\}
\end{align}
%We will provide bounds that compare these against $V_\gn^f$ and $V_\gn^g$, and show that the performance loss approaches zero as the difference between $f$ and $g$ tends to zero, $N\to\infty$ and $\gamma\to1$.
%Note that $\gamma\to1$ will not be required for the cost \eqref{eq:performance1}, which is a discounted cost with the same discount factor that controls are designed for.
%For $f=g$, such suboptimality results were achieved in, e.g., 
%Suboptimality results that quantify the impact of finite horizon and discounting compared to the ideal (usually intractable) infinite-horizon undiscounted case were achieved in, e.g., \cite{grimm_model_2005, tuna2006shorter,grune_analysis_2009,Wort12,grune_economic_2013, postoyan_stability_2017, granzotto_finite-horizon_2021}, and are adapted to the case with plant-model mismatch in this paper.
%For this purpose, we will leverage the stability properties derived in our answer to question 1).

\textit{3) %Continuity/
Robustness of the value function}, i.e.,
how much $V_\gn^g$ and $V_\gn^f$ differ based on the plant-model mismatch $|f-g|_{\mc S}$. %the difference $$ between $f$ and $g$.
%This is relevant because, it is highly desirable to know $V_\gn^g$ for, e.g., policy extraction, performance benchmarking and showing both numerical stability and convergence of algorithms. However, 
Since $g$ is unknown or too complex, calculating $V_\gn^g$ is numerically intractable, and $V_\gn^f$ serves as an approximation to $V_\gn^g$, for which known bounds on the approximation accuracy are important, e.g., for control design.

\begin{table*}[]
    \footnotesize
    \centering
    \hrule
    \begin{align*}
    M_{\gamma,N} &:= \sqrt B\sum_{k=0}^{N-1}\left(\sqrt{\gamma\left(1-\frac1B\right)} L\right)^k
    &K_{\gamma,N} &:= %\min\left\{M_\gn,
    \sum_{k=0}^{N-1}\left(\gamma L^2\right)^k %\right\}
    &F_\No &:= \left(1-\frac1B\right)^{N_0-1}B^2
    \end{align*}
    \vspace{-0.3cm}
    \begin{align*}
    \kappa_\gn(s) &:= K_\gn s^2 + 2M_\gn s
    &\overline\kappa_\gamma(s) &:= \inf_{N_0\in\N}\left(\kappa_\gnn(s) + F_\No(s+1)^2\right)\\
    \overline\kappa_\gn(s) &:= \begin{cases}\min\left\{\kappa_\gn(s), \overline\kappa_\gamma(s)\right\}, &N<\infty\\
    \overline\kappa_\gamma(s),&N=\infty\end{cases}
        &D&:=\left(\sqrt\frac1{\lmin_Q}+\sqrt\frac B{\lmin_R}\right)\sqrt\frac{\lmax_Q\gamma}B\\
        \agn &:= 1- B\okappa_\gn(D|f-g|_{\mc S}) &
        \aggn &:= \frac1\gamma\Big(1- B\okappa_\gn(D|f-g|_{\mc S}) - B(1-\gamma)\Big)\\
        \agnn &:= 1- \underbrace{B\okappa_\gn(D|f-g|_{\mc S})}_{\to0\text{ as }|f-g|_{\mc S}\to0} - \underbrace{B^2e^{-N/B}}_{\to0\text{ as }N\to\infty} &
        \aggnn &:= \frac1\gamma\Big(1- \underbrace{B\okappa_\gn(D|f-g|_{\mc S})}_{\to0\text{ as }|f-g|_{\mc S}\to0} - \underbrace{B(1-\gamma)}_{\to0\text{ as }\gamma\to1} - \underbrace{B^2e^{-N/B}}_{\to0\text{ as }N\to\infty}\Big)
    \end{align*}
    \hrule
    \caption{Definitions of the variables in Sections \ref{s:continuity_of_VF} and \ref{sect:main-results}.}
        \label{tab:variables}\vspace{-0.5cm}
\end{table*}

\section{Continuity of the value function}\label{s:continuity_of_VF}

\begin{figure}
    \centering    
    \begin{flushright}
    \includegraphics[width=\linewidth]{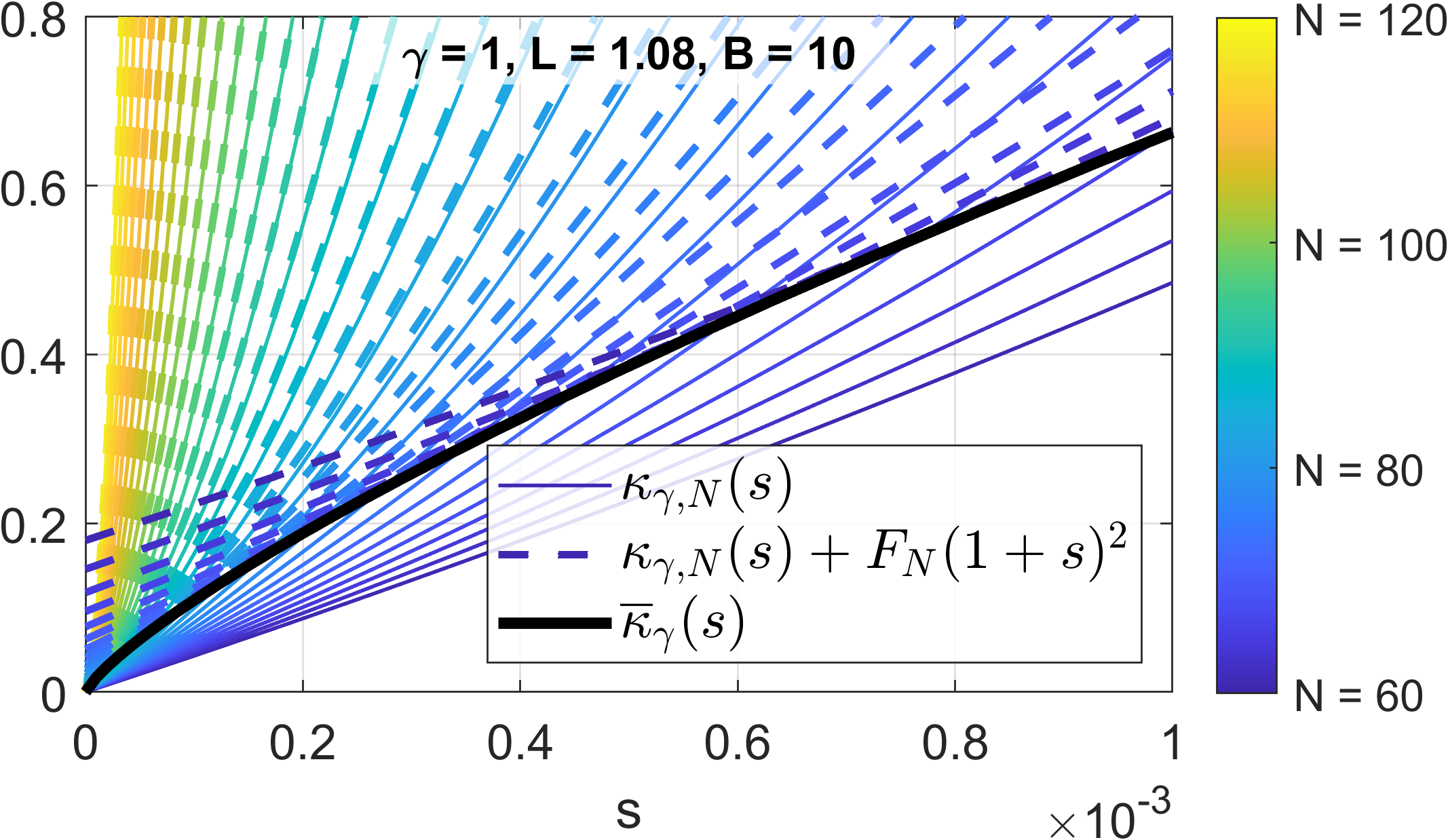}
    \end{flushright}\vspace{-0.2cm}
    \caption{
    Horizon-dependent bounds $\kappa_\gn(s)$ (solid in color) and uniform bound $\okappa_\gamma(s)$ (black), constructed as the infimum (lower envelope) of $\kappa_\gn(s)+F_N(s+1)^2$ (dashed in color) over $N\in\N$.
    %, with $\gamma=1, L=1.08, B=10$.
    }
    \label{fig:visulation}
\end{figure}

Our key tool to %establish %for stability under plant-model mismatch
answer all three questions posed in Section \ref{s:2C} is continuity of the value function, which will later be used as a Lyapunov function.
%Apart from that, continuity is also useful when approximating the value function (via interpolation, neural networks etc.).
We begin with a bound for $V_\gn^f$ that applies only to finite-horizon costs.
\begin{prop}[\textbf{Continuity of value function}] %~$V_{\gn}^N$]
\label{prop:QLipschitz}
    %Suppose that Assumptions~\ref{assm:Lipschitz} and~\ref{assm:cost_controllable} are satisfied.
    %Then, 
    For every %any
    $\gamma\in(0,1]$ and $N\in\N$, there exist constants
    $M_\gn$, $K_\gn$ 
    (given in Table~\ref{tab:variables}) such that, for all $x,y\in\R^n$, 
    \begin{align}
        &V_\gn^f(y) - V_\gn^f(x)\nonumber\\
        &\leq2M_\gn\sqrt{\lmax_Q}\|x\|_Q|x-y|+K_\gn\lmax_Q|x-y|^2.\label{eq:prop:QLipschitz2}
    \end{align}
    If $x\neq0$, this can be formulated with $\kappa_\gn$ in Table \ref{tab:variables} as
    \begin{align}
        V_\gn^f(y) - V_\gn^f(x)  \leq \kappa_\gn\left(\tfrac{\sqrt{\lmax_Q}|x-y|}{\|x\|_Q}\right)\|x\|_Q^2.\label{eq:prop:QLipschitz}
    \end{align}
    %\end{subequations}
    % \vspace{-0.4cm}
    %     \begin{subequations}
    %         \begin{align}
    %             &|V_\gn^f(x) - V_\gn^f(y)|\\
    %             &\leq M_\gn\sqrt{\lmax_Q}|x-y|(\|x\|_Q+\|y\|_Q)\label{eq:prop:QLipschitz_alternative0}\\
    %             &\leq2M_\gn\sqrt{\lmax_Q}\|x\|_Q|x-y|+M_\gn\lmax_Q|x-y|^2\label{eq:prop:QLipschitz_alternative}
    %     \end{align}
    %     \end{subequations}
\end{prop}
The functions $\kappa_\gn$ are illustrated in Fig.~\ref{fig:visulation}.
Inequality~\eqref{eq:prop:QLipschitz2} resembles Lipschitz continuity, but is tailored to the quadratic growth of $V_\gn^f$.
%The terms 
$M_\gn$ and $K_\gn$ bound the accumulation of errors. %along the summated cost.
Discounting is beneficial as it shrinks errors, which reduces $M_\gn$ and $K_\gn$ and smoothens the value function. %leading to smaller $M_\gn$ and $K_\gn$ and a smoother value function.
%The potential benefit of discounting is apparent in $M_\gn$, $K_\gn$ and $\kappa_\gn$.
%In fact, if $\gamma (1-\frac1B)L^2<1$, then $M_\gn$ and $\kappa_\gn(s)$ remain bounded as $N\to\infty$ because discounting and stability counteract the accumulation of prediction errors.
%This generalizes the contractivity in the form of $L<1$ discussed in \cite[Remark 4]{liu_certainty-equivalence_2026} to the milder condition $\gamma (1-\frac1B)L^2<1$.
%Then, a single bound for $V_\gn^f(y)-V_\gn^f(x)$ that applies uniformly over all horizon lengths $N\in\Ni$ (including infinity) is achieved when $M_\gn$ and $K_\gn$ in \eqref{eq:prop:QLipschitz} are replaced by their limits as $N\to\infty$.
In the general case of $\gamma (1-\frac1B)L^2\geq1$, however, we have $\lim_{N\to\infty}M_\gn=\infty$, meaning that $\kappa_\gn$ become infinitely steep as $N\to\infty$, as shown in Fig.~\ref{fig:visulation}.
Still, %even then 
a uniform bound %that is uniform over all 
in~$N \in \Ni$ can be achieved. %as in the following proposition. 
\begin{prop}[\textbf{Continuity of value function uniformly %~$V_\gn^N$ %uniform over all 
in~$N \in \Ni$}]\label{prop:QLipschitz_uniform}
    %Suppose that Assumptions~\ref{assm:Lipschitz} and~\ref{assm:cost_controllable} hold. %are satisfied.
    %Then, f
    For any 
    $\gamma\in(0,1]$, there exists $\overline\kappa_\gamma\in\mc K$ (given in Table \ref{tab:variables}) such that
    \begin{align}
        V_\gn^f(y) - V_\gn^f(x) \leq \overline\kappa_\gamma\left(\tfrac{\sqrt{\lmax_Q}|x-y|}{\|x\|_Q}\right)\|x\|_Q^2 \label{eq:prop:QLipschitz_uniform}
    \end{align}    
    holds for all $N\in\Ni$ and $x,y\in\R^n$ with $x\neq0$.
\end{prop}
Inequality \eqref{eq:prop:QLipschitz_uniform} is achieved by truncating the cost at some arbitrary horizon $N_0\in\N$, applying Proposition \ref{prop:QLipschitz} to $V_\gnn^f$, and bounding the remaining terms in the cost using SA\ref{assm:cost_controllable} with factor $F_\No$ (see Table \ref{tab:variables}).
The function $\okappa_\gamma$ is then the infimum (lower envelope) of these bounds over $N_0\in\N$ as illustrated in Fig.~\ref{fig:visulation}.
%The fact that $N_0\in\N$ can be chosen arbitrarily means that we can take the best out of all the resulting bounds.
%This leads to the infimum over $\No\in\N$ in the definition of $\overline\kappa_\gamma$ as illustrated in Figure~\ref{fig:visulation}.
The non-uniform bound \eqref{eq:prop:QLipschitz} can still be better than \eqref{eq:prop:QLipschitz_uniform} for short horizons. Hence, we take the minimum of the two, and obtain
\begin{align}
    V_\gn^f(y) - V_\gn^f(x)\leq \okappa_\gn\left(\tfrac{\sqrt{\lmax_Q}|x-y|}{\|x\|_Q}\right)\|x\|_Q^2\label{eq:prop:QLipschitz_final}
\end{align}
for $\okappa_\gn$ defined in Table \ref{tab:variables}.
This is the bound used in the remainder of this paper.

\section{Stability, suboptimality and robustness of the value function}\label{sect:main-results}

In this section, we develop the relaxed DP framework and answer the three main questions posed in Section \ref{s:2C}.
%For stability under plant-model mismatch, we aim to establish the value function $V_\gn^f$ as a Lyapunov function for the closed loop ~\eqref{eq:closed_loop}. For this purpose, we derive the next relaxed DP inequalities.
%What is known directly is the Bellman equation \eqref{eq:DP} (see Table \ref{table:rDP}).
%To turn this into a Lyapunov condition, $f(x,u)$ in $\gamma V_\gnm^f(f(x,u))$ must be replaced by $g(x,u)$ levaraging the continuity of the value function, i.e., inequality~\eqref{eq:prop:QLipschitz_final}, $N-1$ be replaced by $N$, and the factor $\gamma$ be removed.
%The following theorem makes (various selections of) these adaptions, as illustrated in Figure \ref{table:rDP}.
\begin{prop}[\textbf{Relaxed Dynamic Programming}]\label{prop:rDP}
    Consider plant \eqref{eq:system-dynamics} and model \eqref{eq:surrogate}. % and suppose that Assumptions~\ref{assm:Lipschitz} and \ref{assm:cost_controllable} hold. %are satisfied.
    Then, for any %Let the 
    discount factor $\gamma\in(0,1]$ and horizon $N\in\Ni, N\geq2$, %be given.
    %Then, 
    the inequalities
    \begin{align}
        \gamma V_\gnm^f(g(x,u)) - V_\gn^f(x) &\leq -\agn\ell(x,u) \tag{RDP1}\label{eq:rDP0} \\
        \gamma V_\gn^f(g(x,u)) - V_\gn^f(x) &\leq -\agnn\ell(x,u) \tag{RDP2}\label{eq:rDPn} \\
        \phantom{\gamma}V_\gnm^f(g(x,u)) - V_\gn^f(x) &\leq -\aggn\ell(x,u) \tag{RDP3}\label{eq:rDPg}\\
        V_\gn^f(g(x,u)) - V_\gn^f(x) &\leq -\aggnn\ell(x,u) \quad \tag{RDP4}\label{eq:rDPgn} 
    \end{align}
    hold for all $x\in\mc S$ and $u\in\mc U_\gn^f(x)$, with \emph{suboptimality indices} $\agn, \agnn, \aggn,\aggnn$ defined in Table \ref{tab:variables}.
\end{prop}

Proposition \ref{prop:rDP} establishes a relaxed DP framework that accounts for plant-model mismatch, horizon length and discounting.
Inequality \eqref{eq:rDPgn} is of main interest for Lyapunov stability.
Inequalities \eqref{eq:rDP0}-\eqref{eq:rDPg} will be used for suboptimality and robustness of the value function.
%The coefficients $\agn, \aggn, \agnn, \aggnn$ etc. are called \emph{suboptimality indices}.
In the ideal case of $f=g$, $N=\infty$ and $\gamma=1$, the suboptimality indices are all equal to 1, and \eqref{eq:rDP0}-\eqref{eq:rDPgn} recover (one side of) \eqref{eq:DP}.
Otherwise, moving from \eqref{eq:DP} to \eqref{eq:rDP0}-\eqref{eq:rDPgn} incurs relaxation terms which reduce the suboptimality indices.
All of \eqref{eq:rDP0}-\eqref{eq:rDPgn} change $f(x,u)$ in \eqref{eq:DP} to $g(x,u)$.
The effect of this is bounded using \eqref{eq:prop:QLipschitz_final} by the relaxation term $B\okappa_\gn(D|f-g|_{\mc S})$.
Since $\okappa_\gn\in\mc K$, we have $\okappa_\gn(D|f-g|_{\mc S})\to0$ as $|f-g|_{\mc S}\to0$.
This means that small plant-model mismatch causes only a small performance loss. %reduction in the suboptimality indices.
The underscores under $N$ and $\gamma$ in the suboptimality indices indicate whether the relaxation terms $B^2e^{-N/B}$ and $B(1-\gamma)$ are present.
%Furthermore, \eqref{eq:rDPn} and \eqref{eq:rDPgn} change $V_\gnm^f$ in \eqref{eq:DP} to $V_\gn^f$, which incurs the relaxation term $B^2 e^{-N/B}$.
%Hence, \eqref{eq:rDPn} and \eqref{eq:rDPgn} feature the suboptimality indices $\agnn$ and $\aggnn$, where the underscore under $N$ indicates that this relaxation term is present.
%Similarly, \eqref{eq:rDPg} and \eqref{eq:rDPgn} remove the factor $\gamma$ and therefore use the suboptimality indices $\aggn$ and $\aggnn$ with the relaxation term $B(1-\gamma)$ linked to discounting.

We now state our main theorem, which answers the three questions posed in Section \ref{s:2C}.
%This is followed by a detailed discussion and corollaries\karl{, which show/showing ...}.
The notation $\mc L_\gn^f(\mc S)$ for $\mc S\subseteq\R^n$ shall denote the largest level set $\{x\in\R^n~|~V_\gn^f(x)\leq c\}, c\geq0$ of $V_\gn^f$ that is contained~in~$\mc S$.
\begin{thm}\label{thm:main}
   Consider plant \eqref{eq:system-dynamics} and model \eqref{eq:surrogate}.
   %and suppose Assumptions~\ref{assm:Lipschitz} and \ref{assm:cost_controllable} are satisfied.
    Then, for any discount factor $\gamma\in(0,1]$ and horizon $N\in\Ni, N\geq2$ satisfying $\aggnn>0$, the following hold.
    \begin{enumerate}
    \item[\textbf{1)}] \textbf{Exponential stability.} %If $\aggnn>0$ or $\mc S=\R^n$, then 
    Any solution $x_k$, $k\in\N_0$, to \eqref{eq:closed_loop} with $x_0\in\mc L_\gn^f(\mc S)$ satisfies
        %\begin{align}
            $\|x_k\|_Q \leq \sqrt B \left(1-\aggnn/B\right)^{k/2}\|x_0\|_Q \label{eq:thm:main:exponential}$
        %\end{align}
        for all $k\in\N_0$.
    \item[\textbf{2)}] \textbf{Suboptimality.} %If $\aggnn\geq0$ or $\mc S=\R^n$, then 
    For any $x\in\mc L_\gn^f(\mc S)$, 
    %the inequalities 
        \begin{align}
        \aggnn \overline J_\oi^g\left(x,\mc U_\gn^f\right) &\leq V_\gn^f(x),\label{eq:thm:suboptimality2}\\
        \agnn \overline J_\gi^g\left(x,\mc U_\gn^f\right) &\leq V_\gn^f(x)\label{eq:thm:suboptimality}
    \end{align}
    hold for the performance indices given in \eqref{eq:performance1} and \eqref{eq:performance2}.
    \item[\textbf{3)}] \textbf{Robustness of the value function.} Let $x\in\R^n$. If one of the conditions
    \begin{enumerate}
        \item[(a)] $\mc S=\R^n$, or
        \item[(b)] $V_\gn^f(x)\leq\overline c$, where $\overline c$ is the largest $c\in\R_{\geq0}\cup\{\infty\}$ with $\{z\in\R^n~|~\|z\|_Q^2\leq c\}\subseteq\mc S$, or
        \item[(c)] $\aggn\!\geq\!\frac{2B\log(B)}{N}$, $B\geq3$ and $x\in\mc L_\gn^f(\mc S)$
    \end{enumerate}
    holds, then
    %\begin{align}
        $\agn V_\gn^g(x) \leq V_\gn^f(x). $%\label{eq:thm:robustness_value_function}
    %\end{align}
    \end{enumerate}
\end{thm}
% Theorem \ref{thm:main} ensures %guarantees 
% horizon length %$N$ 
% and plant-model mismatch. %$|f-g|_{\mc S}$.
% %This dependency is expressed through the suboptimality indices.

Theorem \ref{thm:main} guarantees closed-loop stability in item 1), bounds for closed-loop performance in item 2), and bounds between the value functions for $f$ and $g$ in item 3), all depending on the discount factor $\gamma$, horizon length $N$ and plant-model mismatch $|f-g|_{\mc S}$.

All suboptimality indices include the term $\okappa_\gn(D|f-g|_{\mc S})$, which bounds the effect of plant-model mismatch. Since $\okappa_\gn\leq\okappa_\gamma$ and $\okappa_\gamma$ is independent of $N$, our results are uniform in $N$ in the sense that a longer horizon does not require increasingly smaller plant-model mismatch.
This is important when a long horizon is required, e.g., %/required 
in MPC.
To the best of our knowledge, this is novel in this context.
In contrast, existing results in \cite{SchiWort25, liu_certainty-equivalence_2026} are horizon-dependent, with bounds that generally deteriorate for longer horizon and fail over infinite horizon.
%We demonstrate this on an example in Section \ref{sect:illustrative-example} below.
In addition, our results generalize those in \cite{SchiWort25, liu_certainty-equivalence_2026} by considering discounted costs.
We demonstrate the benefits of the uniformity over $N$ as well as discounting in the example in Section \ref{sect:illustrative-example} below.

The following corollary states the immediate consequence of item 3) of Theorem \ref{thm:main} that cost-controllability of the model~\eqref{eq:surrogate} implies cost-controllability of the plant~\eqref{eq:system-dynamics} when plant-model mismatch is small enough.
\begin{corollary}[\textbf{Cost-controllability under plant-model mismatch}]\label{cor:cost_controllability_plant}
    Suppose that %Assumptions \ref{assm:Lipschitz} and \ref{assm:cost_controllable} hold, and that 
    $\mc S=\R^n$ and $\alpha_\oi>0$. Then, for any $x\in\R^n$, it holds that
    %\begin{align}
        $V_\oi^g(x) \leq \nicefrac{B}{\alpha_\oi} \|x\|_Q^2.$
    %\end{align}
\end{corollary}

The following corollary provides alternative formulations of \eqref{eq:thm:suboptimality2} and \eqref{eq:thm:suboptimality} that compare against $V_\gn^g$ instead of $V_\gn^f$.
We obtain this by combining item 2) of Theorem~\ref{thm:main} with item 3) of Theorem \ref{thm:main} with $f$ and $g$ swapped.
Note that cost-controllability of $g$ with %$B\alpha_\oi^{-1}$ 
$\nicefrac{B}{\alpha_\oi}$ is ensured by Corollary \ref{cor:cost_controllability_plant}.
For simplicity, we formulate this for the case $\mc S=\R^n$.
\begin{corollary}[\textbf{Suboptimality vs. $g$}]\label{cor:suboptimality_vsg}
    Suppose that %Assumptions \ref{assm:Lipschitz} and \ref{assm:cost_controllable} hold,
    $\mc S=\R^n$, $\alpha_\oi>0$, $g$ is continuous, and $|g(x,u)-g(y,u)|\leq L|x-y|$ for all $x,y\in\R^n$ and $u\in\mathbb U$.
    Further, let $\gamma\in(0,1]$ and $N\in\Ni, N\geq2$ with $\aggnn>0$. Let $\widetilde\alpha_\gn$ be defined as $\agn$ but with $\nicefrac{B}{\alpha_\oi}$ instead of $B$ (including in $\okappa_\gn$ and terms therein).
    Then, for every $x\in\R^n$,
    \begin{align}
        \widetilde\alpha_\gn  \aggnn \overline J_\oi^g\big(x,\mc U_\gn^f\big) &\leq V_\gn^g(x),\label{eq:thm:suboptimality2_vsg}\\
        \widetilde\alpha_\gn \agnn \overline J_\gi^g\big(x,\mc U_\gn^f\big) &\leq V_\gn^g(x).\label{eq:thm:suboptimality_vsg}
    \end{align}
\end{corollary}

\begin{figure}
    \centering
    \includegraphics[width=0.97\linewidth]{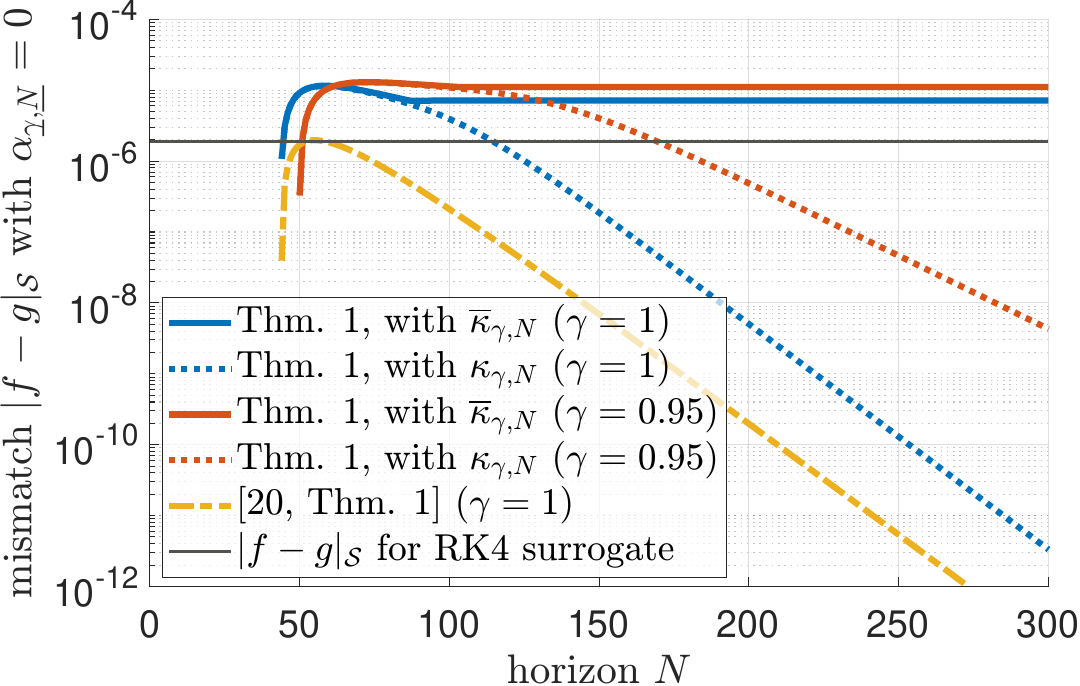}
    \caption{Comparison of bounds for the example in Section \ref{sect:illustrative-example}.
    }
    \label{fig:example}
\end{figure}

\section{Illustrative example}\label{sect:illustrative-example}
%--
%\rc{Title: Illustrative examples or case studies?}
The purpose of this example is to demonstrate the bounds of Theorem~\ref{thm:main} and compare them against \cite{liu_certainty-equivalence_2026}. %, and demonstrate the uniformity over $N$.
Consider an inverted pendulum with nonlinear continuous-time dynamics
\begin{align}
    \dot x_1=x_2,\quad \dot x_2=1.25\sin(x_1)-x_2+u, \label{eq:exmp:nonlinear}
\end{align}
where $x_1$ is the angle between the rod and the vertical axis  ($x_1=0$ corresponding to the upward position), $x_2$ is the angular velocity, and $u$ is the control input.
%, $\mathfrak g=1$ is the ratio between gravitational acceleration and length of the rod, and $\mathfrak d=1$ is a damping coefficient.
Let the plant \eqref{eq:system-dynamics} be the exact discretization of \eqref{eq:exmp:nonlinear} via zero-order hold input (which cannot be derived analytically) and let the surrogate model \eqref{eq:surrogate} be its 4th-order Runge-Kutta discretization, both with step size $T=0.1$.
Since $g$ is unknown, $f$ is used to solve the optimal control problem.
We choose $\mathbb U=[-20,20]$, $\mc S=[-\pi/2,\pi/2]\times[-2,2]$, $Q=\begin{bsmallmatrix}10&0\\0&1\end{bsmallmatrix}, R=0.1$, %i.e., deviations of up to $90^\circ$ from the upright position,
which satisfy SA\ref{assm:1:setU} and SA\ref{assm:4:QR}.
SA\ref{assm:Lipschitz} (Lipschitz continuity) is satisfied with $L\approx1.076$ calculated as the largest spectral norm of the Jacobian $\partial f/\partial x$ on $\mc S\times\mathbb U$.
For SA\ref{assm:cost_controllable} (cost-controllability), we simulated the surrogate under a saturated LQR feedback on a grid of initial states in $\mc S$, which yields $B\approx 9.616$.
To verify SA\ref{assm:3:mismatch} (finite proportional mismatch) and determine $|f-g|_{\mc S}$, we approximated the true unknown $g$ to machine precision with the 4th-order Runge-Kutta method over 1000 subdivisions of the discretization interval, and then found $|f-g|_{\mc S}\approx1.5\cdot10^{-6}$ by evaluating \eqref{eq:defn:mismatch} on a fine grid.

Figure \ref{fig:example} shows, as a function of $N$, the threshold on the mismatch $|f-g|_{\mc S}$ below which stability is guaranteed by Theorem \ref{thm:main}, i.e., the value at which $\aggnn=0$.
%Figure \ref{fig:example} shows, as a function of $N$, the largest mismatch $|f-g|_{\mc S}$ for which stability is still guaranteed, i.e., the value at which  $\aggnn=0$. %Figure \ref{fig:example} shows the stability threshold, i.e., the mismatch $|f-g|_{\mc S}$ for which $\aggnn=0$, as function of $N$.
%By Theorem \ref{thm:main}, stability is guaranteed for any mismatch below this threshold.
Note that this threshold can only exist for $N\geq2B\log(B)$ because otherwise $\aggnn<0$ even if $f=g$. %for small $N$, $\aggnn<0$ even if $f=g$.
We compare the uniform bound for $\aggnn$ defined using $\okappa_\gn$ (solid lines) versus the non-uniform bound that would be obtained if $\aggnn$ were defined with $\kappa_\gn$ instead (dotted lines). The latter, as well as the bound in \cite{liu_certainty-equivalence_2026}\footnote{We used, in the notation of \cite{liu_certainty-equivalence_2026}: $m_{\ell,x}=2$, $L_{\ell,x}=20$, $m_{\ell,u}=L_{\ell,u}=0.2$, $M_{\mathrm{ec}}=11$, $\bar\gamma=B-1\approx8.62$, $\varepsilon_F\approx1.26$, $L_{f,x}=L\approx1.076$, $L_{f,u}\approx0.096$, $L_\Delta(\delta(\theta))=|f-g|_{\mc S}$, $H\approx10$, $\lambda=0.77$.} %, the last two being related to the assumptions in \cite{liu_certainty-equivalence_2026} mentioned in Remark \ref{rem:comparison_liu}} %. (from the LQ problem of the linearized surrogate, chosen to maximize the admissible mismatch).}
(yellow line) deteriorate as $N$ increases.
Consequently, our bound guarantees stability for the mismatch $|f-g|_{\mc S}\approx1.5\cdot10^{-6}$ found earlier (black line) for all horizons $N\geq45$, whereas \cite[Theorem 1]{liu_certainty-equivalence_2026} guarantees stability only for a small range of~$N$.
While our bounds may still be fairly conservative, we plan to reduce conservatism by exploiting error contraction under a local LQR approximation in future research. %, which can lead to numerical guarantees that are practically useful in a wide range of scenarios, all still within the relaxed DP framework presented in this paper.

The example also demonstrates the benefit of discounting by comparing bounds for $\gamma=1$ (blue in Fig. \ref{fig:example}) versus $\gamma=0.95$ (red in Fig. \ref{fig:example}).
While discounting with $\gamma=0.95$ consumes almost half of the robustness margin through the term $B(1-\gamma)\approx0.481$ in $\aggnn$, it reduces $\okappa_\gn$ enough to provide a net benefit for robustness, as the uniform stability threshold for $\gamma=0.95$ (red horizontal line) lies above the undiscounted one (blue horizontal line).

% Fig. \ref{fig:example} shows the resulting suboptimality indices.
% With the horizon-specific bound $\kappa_{1,N}$ (red curve) the suboptimality index increases at first, but eventually decreases and becomes negative because $M_\gn\to\infty$ as $N\to\infty$.
% For $\alpha_{\underline 1,\underline N}$ (which is the same as $\alpha_{1,\underline N}$ defined in Table \ref{tab:variables} (yellow curve) this is not the case, and stability and suboptimality are guaranteed for arbitrarily long horizon (without the need to use more and more accurate models).
% Finally, we compare with the suboptimality index in \cite{liu_certainty-equivalence_2026}, which also falls of for large $N$ and therefore can guarantee stability only for a range of $N$.
% Furthermore, our bound improves upon \cite{liu_certainty-equivalence_2026} in the sense that our suboptimality indices feature the term $B^2e^{-N/B}$ which is exponentially decaying in $N$, whereas the corresponding term in \cite{liu_certainty-equivalence_2026} is of order $\mathcal{O}(B^2/N)$ and, thus, decays only linearly.
% The reduction in the required horizon length for stability allows significantly larger plant-model mismatch compared to \cite{liu_certainty-equivalence_2026}, as seen in Fig. \ref{fig:example}.

%--
\section{Conclusion}\label{sect:conclusion}
%--
We have presented stability and suboptimality guarantees for plants controlled by optimal inputs generated using a surrogate model. 
The main novelties are the combined relaxed DP framework, the uniformity of the bounds over all (possibly infinite) horizons, and the benefit of discounting to handle plant-model mismatch.
The results may be a foundation for future research where optimal controls are directly learned from data (i.e., reinforcement learning).
Then, the relaxed DP framework may be expanded by parameters (in addition to, or replacing $N, \gamma$ and $|f-g|_{\mc S}$) such as the amount and quality of data, number of training steps, learning rate, and richness of an approximation architecture.
The uniformity of our results in $N$ may then correspond to uniformity in the number of training steps (e.g., steps of value iteration), assuring that stability counteracts error accumulation. % (with or without discounting).

\bibliographystyle{ieeetr}
\bibliography{references}

\appendix

%\section{Proofs of Section \ref{s:continuity_of_VF}}\label{s:proof:prop:QLipschitz}
%\label{s:proof:prop:all_horizons}
\noindent\textbf{Proof of Proposition \ref{prop:QLipschitz}.}
    %Given $L,B\geq0$ and $N\in\N$, we choose
    Let %$f,L$ and $B$ be as in the statement of Proposition \ref{prop:QLipschitz}, and let %$f\in\mc M_{L,B}^0,
    $N\in\N$, $\gamma\in(0,1]$ and $x,y\in\R^n$ be given.
    Proposition \ref{prop:1.bellman} implies the existence of $\mathbf u_N\in\mathbb U^N$ satisfying $V_\gn^f(x) = J_\gn^f(x,\mathbf u_N)$. 
    Define $x_k:=\varphi^f_x(k,\mathbf u_{k})$ and $y_k:=\varphi^f_y(k,\mathbf u_{k})$ for $k\in\{0,\dots,N-1\}$.
    Then, by \eqref{eq:L-Lipschitz},
    \begin{align}\label{eq:proof:thm:QLipschitz_1}
        |x_k-y_k| \leq %L^k|x_0-y_0| = 
        L^k|x-y|\quad\forall  k\in\{0,\dots,N-1\}.
    \end{align}
    Further, for all $k\in\{0,\dots,N-2\}$, by choice of $\mathbf u_N$, Proposition \ref{prop:1.bellman} and \eqref{eq:B-cost-controllable},
    \begin{align}
        \gamma V_{\gamma,N-k-1}^f(x_{k+1}) &= V_{\gamma,N-k}^f(x_k) -\ell(x_k,u_k)\nonumber\\
        &\leq \left(1-\tfrac1B\right) V^f_{\gamma,N-k}(x_k).
    \end{align}
    %and $\gamma^k \ell(x_k,u_k) \leq J_{\gamma,N}^f(x,u_{0:N-1})$ yield
    % \begin{align}
    %     |x_k|^2 &\leq \frac1{\underline\lambda}\|x_k\|_Q^2 \leq \frac1{\underline\lambda} \ell(x_k,u_k)\nonumber\\
    %     &\leq \frac1{\underline\lambda \gamma^{k}}V_\gn^f(x)
    %     \leq \frac B{\underline\lambda \gamma^k}\|x\|_Q^2
    %     \leq \frac{\overline\lambda B}{\underline\lambda \gamma^k}|x|^2 \label{eq:proof:thm:QLipschitz_2}
    % \end{align}
    Iterating this, we get for all $k\in\{0,\dots,N-1\}$ that, by~\eqref{eq:B-cost-controllable},
    \begin{align}
        \|x_k\|_Q^2 \leq V^f_{\gamma,N-k}(x_k)&\leq \left(\tfrac1\gamma\left(1-\tfrac1B\right)\right)^k V_\gn^f(x)\nonumber\\
        &\leq \left(\tfrac1\gamma\left(1-\tfrac1B\right)\right)^k B\|x\|_Q^2.
        \label{eq:proof:thm:QLipschitz_2}
    \end{align}
    Then, using $V_\gn^f(x) = J_\gn^f(x,\mathbf u_N)$, \eqref{eq:proof:thm:QLipschitz_1} and \eqref{eq:proof:thm:QLipschitz_2},
    and writing $\oll := \lmax_Q$ for brevity,
    % On the other hand, using symmetry and positive definiteness of $Q$, we get 
    % \begin{align}
    %     &|\ell(x,u) - \ell(y,u)| = |x^\top Qx - y^\top Qy| = |(x-y)^\top Q(x+y)|\nonumber\\
    %     &\leq\overline\lambda|x-y\|x+y|\leq\oll|x-y|(|x|+|y|). \label{eq:proof:prop:QLipschitz_binomial}
    % \end{align}
    % Using the definition of $V_{\gamma,N}(y)$ and combining inequalities~\eqref{eq:proof:thm:QLipschitz_1}, \eqref{eq:proof:thm:QLipschitz_2} and \eqref{eq:proof:prop:QLipschitz_binomial} yields
    % \begin{align*}
    %     & V_\gn(x) - V_\gn(y) \leq J_\gn(x,u_{0:N-1}) - J_\gn(y,u_{0:N-1})\\
    %     = & \sum\nolimits_{k=0}^{N-1}\gamma^k \underbrace{( \ell(x_k,u_k)-\ell(y_k,u_k))}_{\leq \overline\lambda|x_k-y_k|(|x_k|+|y_k|) \leq \overline\lambda|x_k-y_k| (2|x_k| + |x_k-y_k|)} \\
    %     \leq & \sum\nolimits_{k=0}^{N-1}\gamma^k \overline{\lambda} L^k |x-y| \Big( 2\sqrt{ \frac {\overline{\lambda} B}{\underline{\lambda} \gamma^{k}}} |x| + L^k \underbrace{|x-y|}_{\leq |x|+|y|} \Big) \\
    %     \leq & |x-y| \Big( |x|+|y| \Big) \cdot \sum\nolimits_{k=0}^{N-1}\gamma^k\overline\lambda L^k \Big( 2\sqrt{\frac{\overline\lambda B}{\underline\lambda \gamma^k}} + L^k \Big) \\
    %     \stackrel{\gamma \leq 1}{\leq} & M_N|x-y|(|x|+|y|).\numberthis
    % \end{align*}
    \begin{align}
        &\quad~ V_\gn^f(y) - V_\gn^f(x) 
        \leq J_\gn^f(y,\mathbf u_N) - J_\gn^f(x,\mathbf u_N)\nonumber\\
        &=\sum_{k=0}^{N-1}\gamma^k(\ell(y_k,u_k) - \ell(x_k,u_k))\nonumber\\
        %&=\sum_{k=0}^{N-1}\gamma^k(y_k^\top Qy_k - x_k^\top Qx_k)\nonumber\\
        %&=\sum_{k=0}^{N-1}\gamma^k(2x_k+y_k-x_k)Q(y_k-x_k)\nonumber\\
        &=2\sum_{k=0}^{N-1}\gamma^k x_k^\top Q(y_k-x_k) + \sum_{k=0}^{N-1}\gamma^k (y_k-x_k)^\top Q(y_k-x_k)\nonumber\\
        &\leq2\sum_{k=0}^{N-1}\gamma^k\sqrt\oll \|x_k\|_Q|x_k-y_k| + \sum_{k=0}^{N-1}\gamma^k\oll|x_k-y_k|^2\nonumber\\
        &\leq 2\sqrt\oll\sum_{k=0}^{N-1}\gamma^k\left(\tfrac1\gamma\left(1-\tfrac1B\right)\right)^{k/2}\sqrt B\|x\|_Q L^k|x-y|\nonumber\\
        &~~~~+ \oll\sum_{k=0}^{N-1}\gamma^kL^{2k}|x-y|^2\nonumber\\
        &= 2M_{\gn}\sqrt\oll \|x\|_Q|x-y| + K_{\gn}\oll|x-y|^2,\label{eq:prop:proof:QLipschitz}
    \end{align}
    completing the proof of \eqref{eq:prop:QLipschitz2}, from which \eqref{eq:prop:QLipschitz} follows. 
    %Inequality \eqref{eq:prop:QLipschitz} follows immediately given the definition of $\kappa_\gn$ in Table~\ref{tab:variables}.
    % To prove \eqref{eq:prop:QLipschitz_alternative}, define for $m\in\N$ the states $x_i := x + \frac im(y-x)$.
    % By triangle inequality, $\|x_i\|_Q\leq \frac{m-i}m\|x\|_Q + \frac{i}m\|y\|_Q$.
    % Using this and \eqref{eq:prop:proof:QLipschitz},
    % \begin{align}
    %     &V_\gn^f(y) - V_\gn^f(x)
    %     = \sum_{i=0}^{m-1} \left(V_\gn^f(x_{i+1})-V_\gn^f(x_i)\right)\nonumber\\
    %     &\leq \sum_{i=0}^{m-1}\left(2M_\gn\sqrt{\oll}\|x_i\|_Q|x_{i}-x_{i+1}| + K_\gn\oll|x_{i}-x_{i+1}|^2\right)\nonumber\\
    %     &= \sum_{i=0}^{m-1}\left(2M_\gn\sqrt{\oll}\|x_i\|_Q\frac{|x-y|}m + K_\gn\oll\frac{|x-y|^2}{m^2}\right)\nonumber\\
    %     &\leq 2M_\gn\sqrt\oll\frac{|x-y|}m\left(\frac{m(m+1)}{2m}\|x\|_Q + \frac{m(m-1)}{2m}\|y\|_Q\right) \nonumber \\
    %     &~~~~+ K_\gn\oll\frac{|x-y|^2}m.
    % \end{align}
    % Taking the limit as $m\to\infty$ yields \eqref{eq:prop:QLipschitz_alternative}.

\noindent\textbf{Proof of Proposition \ref{prop:QLipschitz_uniform}.}
Let $N\in\Ni$, $\gamma\in(0,1]$ and $x,y\in\R^n$ with $x\neq0$ be given.
Furthermore, let $N_0\in\N$ be arbitrary.
Proposition~\ref{prop:1.bellman} implies the existence of $\mathbf u_\No\in\mathbb U^\No$ satisfying $V_\gnn^f(y) = J_\gnn^f(y,\mathbf u_\No)$.
First consider the case where $N\geq\No$. Using Proposition \ref{prop:1.bellman}, followed by \eqref{eq:B-cost-controllable} and \eqref{eq:proof:thm:QLipschitz_2} for $k=N_0-1$, we obtain
\begin{align}
    V_\gn^f(y)&\leq J_{\gamma,N_0-1}^f(y, \mathbf u_{\No-1})\nonumber\\
    &~~~~+ \gamma^{\No-1} V_{\gamma,N-(N_0-1)}^f(\varphi^f_y(N_0-1,\mathbf u_{\No-1})) \nonumber\\
    &\hspace{-0.5cm}\leq J_\gnn^f(y,\mathbf u_\No)+\gamma^{\No-1} B\|\varphi^f_y(N_0-1,\mathbf u_{\No-1})\|_Q^2\nonumber\\
    &\hspace{-0.5cm}\leq  V_\gnn^f(y) + F_\No\|y\|_Q^2. \label{eq:proof:prop:uniform}
\end{align}
%with $F_\No$ defined in Table \ref{tab:variables}.
%If $N<N_0$, then $V_\gn^f(y)\leq V_\gnn^f(y)+F_\No\|y\|_Q^2$  holds as well.
\vspace{-0.1cm}
Define for brevity $s:=\sqrt{\lmax_Q}|x-y| / \|x\|_Q$, then, by triangle inequality, $\|y\|_Q^2 \leq (\|x\|_Q + \sqrt{\lmax_Q}|x-y|)^2 = (1+s)^2\|x\|_Q^2$.
Combining this with \eqref{eq:proof:prop:uniform}, $N\geq N_0$ and Proposition \ref{prop:QLipschitz}, %if $N\geq N_0$,
\begin{align}
    &V_\gn^f(y) - V_\gn^f(x)
    \leq V_\gnn^f(y) + F_\No\|y\|_Q^2 - V_\gnn^f(x)\nonumber\\
   &\leq \left(\kappa_\gnn(s) + F_\No(1+s)^2\right)\|x\|_Q^2.\label{eq:proof:prop:uniform:2}
\end{align}
%If $N<N_0$, then, by Proposition \ref{prop:QLipschitz}, $V_\gn^f(y)-V_\gn^f(x) \leq \kappa_\gn(s)\|x\|_Q^2 \leq \kappa_\gnn(s)\|x\|_Q^2$, hence the bound in \eqref{eq:proof:prop:uniform:2} also applies in this case.
If $N<N_0$, then Proposition \ref{prop:QLipschitz} and the fact that $K_\gn$ and $M_\gn$ are monotone in $N$ yield $V_\gn^f(y)-V_\gn^f(x)\leq \kappa_\gn(s)\|x\|_Q^2\leq \left(\kappa_\gnn(s) + F_\No(1+s)^2\right)\|x\|_Q^2$, hence the upper bound in \eqref{eq:proof:prop:uniform:2} holds in this case as well.
Since $N_0\in\N$ is arbitrary, $V_\gn^f(y)-V_\gn^f(x)$ is upper bounded by the infimum of the right-hand sides of \eqref{eq:proof:prop:uniform:2} over $N_0\in\N$, which implies \eqref{eq:prop:QLipschitz_uniform} by definition of $\okappa_\gamma$.
%The proof that $\okappa_\gamma\in\mc K$ can be found in \cite{moldenhauer2026discounted}.

We proceed to prove that $\okappa_\gamma\in\mc K$.
Let $h_\No(s):=\kappa_\gnn(s)+F_\No(s+1)^2$ for $N_0\in\N$ and $s\geq0$, so that $\okappa_\gamma=\inf_{N_0\in\N}h_\No$.
By Table~\ref{tab:variables}, $K_\gnn\geq1$, $M_\gnn\geq\sqrt B$, $0\leq F_\No\leq B^2$, and $F_\No\to0$ as $N_0\to\infty$. Moreover, for all $0\leq t\leq s$,
\begin{align}
    &h_\No(s)-h_\No(t)\label{eq:proof:K:1}\\
    &= (s-t)\big(K_\gnn(s+t)+2M_\gnn+F_\No(s+t+2)\big).\nonumber
\end{align}

\emph{Zero at zero:} $\okappa_\gamma(0)=\inf_{N_0\in\N}F_\No=0$.

\emph{Strictly increasing:} For $0\leq t<s$, \eqref{eq:proof:K:1} gives $h_\No(s)\geq h_\No(t)+2\sqrt B(s-t)\geq\okappa_\gamma(t)+2\sqrt B(s-t)$ for all $N_0\in\N$, so $\okappa_\gamma(s)\geq\okappa_\gamma(t)+2\sqrt B(s-t)>\okappa_\gamma(t)$.

\emph{Continuity:} At $s=0$, $0\leq\okappa_\gamma(s)\leq h_\No(s)\to F_\No$ as $s\searrow0$ for every $N_0\in\N$, and $\inf_{N_0\in\N}F_\No=0$.
On $(0,\infty)$, we show that $\okappa_\gamma$ is Lipschitz on every $[a,b]$ with $0<a<b$.
Let $a\leq t<s\leq b$, $\varepsilon\in(0,1]$, and pick $N_0$ with $h_\No(t)\leq\okappa_\gamma(t)+\varepsilon$. Then $K_\gnn a^2+2M_\gnn a\leq h_\No(t)\leq c:=\okappa_\gamma(b)+1$, so $K_\gnn\leq c/a^2$ and $M_\gnn\leq c/(2a)$. With \eqref{eq:proof:K:1}, $\okappa_\gamma(s)-\okappa_\gamma(t)\leq h_\No(s)-h_\No(t)+\varepsilon\leq C(s-t)+\varepsilon,$
where $C:=2bc/a^2+c/a+2B^2(b+1)$ does not depend on $N_0$, $t$, $s$, or $\varepsilon$. Letting $\varepsilon\to0$ yields $0<\okappa_\gamma(s)-\okappa_\gamma(t)\leq C(s-t)$.
Hence, $\okappa_\gamma$ is continuous.

%\subsection{Proof of Theorem \ref{thm:main}}
%\label{a:proof:main}

\noindent\textbf{Proof of Proposition \ref{prop:rDP}.} Let $\gamma\in(0,1]$, $N\in\N\cup\{\infty\}$ with $N\geq 2$, $x\in\mc S$ and $u\in\mc U_\gn^f(x)$.
    We write for brevity $x^+:=f(x,u)$ and $d:=g(x,u)-f(x,u)$.
    By \eqref{eq:B-cost-controllable} and $u\in\mc U_\gn^f(x)$,
    \begin{align}
        \lmin_R|u|^2 \mkern-1mu\leq\mkern-1mu \|u\|_R^2 \mkern-1mu\leq\mkern-1mu \ell(x,u) \mkern-1mu\leq\mkern-1mu V_\gn^f(x)
        \mkern-1mu\leq\mkern-1mu B \|x\|_Q^2.\label{eq:proof:thm:main:0}
    \end{align}
    Define $\zeta := \sqrt\frac1{\lmin_Q} + \sqrt{\frac B{\lmin_R}}$. Then, using $\|x\|_Q\geq\sqrt{\lmin_Q}|x|$ along with \eqref{eq:proof:thm:main:0} and finally \eqref{eq:defn:mismatch} and $x\in\mc S$,
    \begin{align}
        |d| &= \tfrac{(|x|+|u|)|d|}{|x|+|u|}
        \leq \tfrac{\zeta\|x\|_Q|d|}{|x|+|u|}
        \leq \zeta|f-g|_{\mc S}\|x\|_Q.
        \label{eq:proof:thm:main:1}
    \end{align}
    Furthermore, because $u\in\mc U_\gn^f(x)$ and $N\geq2$, we have $V_\gn^f(x) = \ell(x,u) + \gamma V_\gnm^f(x^+) \geq \ell(x,u) + \gamma\|x^+\|_Q^2 \geq \gamma\|x^+\|_Q^2$, and therefore, with  \eqref{eq:B-cost-controllable},
    \begin{align}
        \|x^+\|_Q^2
        \leq \tfrac{1}{\gamma} V_\gn^f(x)
        \leq \tfrac{B}{\gamma}\|x\|_Q^2.\label{eq:proof:thm:main:2}
    \end{align}
    Overall, using Proposition \ref{prop:QLipschitz_uniform}, followed by the definition of $\okappa_\gamma$ in Table \ref{tab:variables}, \eqref{eq:proof:thm:main:1} and \eqref{eq:proof:thm:main:2},
    and writing $\oll := \lmax_Q$,
    \begin{align*}
        &\gamma V_\gn^f(g(x,u)) - \gamma V_{\gamma,N}^f(x^+)
        \leq \gamma \okappa_\gamma\left(\tfrac{\sqrt{\oll}|d|}{\|x^+\|_Q}\right)\|x^+\|_Q^2\\
        &\leq \gamma \inf_{N_0\in\N}\Big( (K_\gnn+F_\No)\oll|d|^2\\
        &~~~~+ 2(M_\gnn + F_\No)\sqrt{\oll}\|x^+\|_Q|d| + F_\No\|x^+\|_Q^2 \Big)\\
        &\leq \gamma \inf_{N_0\in\N}\Big( (K_\gnn+F_\No)\oll\zeta^2|f-g|_{\mc S}^2\\
        &~~~~+ 2(M_\gnn + F_\No)\sqrt{\oll}\sqrt{\tfrac{B}{\gamma}} \zeta|f-g|_{\mc S} + F_\No\tfrac B\gamma\Big)\|x\|_Q^2\\
        &=B \okappa_\gamma\left(\zeta\sqrt{\tfrac{\oll\gamma}{B}}|f-g|_{\mc S}\right)\|x\|_Q^2.\numberthis\label{eq:proof:rDP1}
    \end{align*}
    If $N\in\N$, we can similarly show using Proposition \ref{prop:QLipschitz} that
    \begin{align}
        &\gamma V_\gn^f(g(x,u)) - \gamma V_{\gamma,N}^f(x^+)\nonumber\\
        &\leq B \kappa_\gn\left(\zeta\sqrt{\tfrac{\oll\gamma}{B}}|f-g|_{\mc S}\right)\|x\|_Q^2.\label{eq:proof:rDP2}
    \end{align}
    Combining \eqref{eq:proof:rDP1} and \eqref{eq:proof:rDP2},
    \begin{align*}
        &\gamma V_\gn^f(g(x,u)) - \gamma V_{\gamma,N}^f(x^+)\\
        &\leq B \okappa_\gn\left(\zeta\sqrt{\tfrac{\oll\gamma}{B}}|f-g|_{\mc S}\right)\|x\|_Q^2\\
        &\leq B\okappa_\gn(D|f-g|_{\mc S})\ell(x,u).\numberthis \label{eq:lem:1}
    \end{align*}
    Since $u\in\mc U_\gn^f(x)$, Proposition \ref{prop:1.bellman} implies $\gamma V_\gnm^f(x^+) - V_\gn^f(x) = -\ell(x,u)$.
    Combining this with \eqref{eq:lem:1} for $N-1$ instead of $N$ and applying $\okappa_\gnm\leq\okappa_\gn$ yields \eqref{eq:rDP0}.
    %To prove \eqref{eq:rDPn}, we aim to upper bound $\gamma V_\gn^f(x^+) - V_\gn^f(x)$. 
    If $N=\infty$, then \eqref{eq:rDPn} is identical to \eqref{eq:rDP0}.
    We proceed to prove \eqref{eq:rDPn} if $N<\infty$.
    Because $u\in\mc U_\gn^f(x)$ and by Proposition \ref{prop:1.bellman}, there exist $u_0=u$ and $u_1,\dots,u_{N-1}\in\mathbb U$ such that, for $\mathbf u_N=(u_0,\dots,u_{N-1})$, $V_\gn^f(x) = J_\gn^f(x,\mathbf u_N)$ holds.
    Denote $x_k:=\varphi^f_x(k,\mathbf u_k)$ for $k\in\{0,\dots,N-1\}$.
    Then, by Proposition \ref{prop:1.bellman},
    \begin{align}
        &\gamma V_\gn^f(x^+) - V_\gn^f(x)\leq\nonumber\\
        &\sum_{k=1}^{N-2}\gamma^k\ell(x_k,u_k) + \gamma^{N-1}V_{\gamma,2}^f(x_{N-1}) - \sum_{k=0}^{N-1}\gamma^k\ell(x_k,u_k)\nonumber\\
        &\leq -\ell(x,u) + \gamma^{N-1}B\|x_{N-1}\|_Q^2 - \gamma^{N-1}\ell(x_{N-1},u_{N-1})\nonumber\\
        &\leq -\ell(x,u) + \gamma^{N-1}(B-1)\|x_{N-1}\|_Q^2. \label{eq:proof:prop4_2}
    \end{align}
    Using \eqref{eq:proof:thm:QLipschitz_2} to bound $\|x_{N-1}\|_Q^2$, as well as $1-\tfrac1B \leq e^{-1/B}$,
    \begin{align}
        &\gamma^{N-1}(B-1)\|x_{N-1}\|_Q^2 \leq \left(1-\tfrac1B\right)^{N-1}B(B-1)\|x\|_Q^2\nonumber\\
        &= B^2\left(1-\tfrac1B\right)^N\|x\|_Q^2 \leq B^2e^{-N/B}\|x\|_Q^2.\label{eq:proof:prop4_3}
    \end{align}

    Adding \eqref{eq:lem:1} and \eqref{eq:proof:prop4_2} and combining this with \eqref{eq:proof:prop4_3} and $\ell(x,u)\geq\|x\|_Q^2$
    %Combining \eqref{eq:lem:1} with $\ell(x,u)\geq\|x\|_Q^2$ and \eqref{eq:lem:2} 
    yields inequality \eqref{eq:rDPn}.    
    Finally, adding $(1-\gamma)V_\gn^f(x)$ to \eqref{eq:rDP0} and then dividing by $\gamma$ yields
    \begin{align*}
        &V_\gnm^f(g(x,u)) - V_\gn^f(x)\\
        &\leq \tfrac1\gamma\left(-\agn\ell(x,u) + (1-\gamma)V_\gn^f(x)\right)\\
        &\leq \tfrac1\gamma\left(-\agn + B(1-\gamma)\right)\ell(x,u) = -\aggn\ell(x,u),\numberthis
        %V_\gn^f(g(x,u)) - V_\gn^f(x) \leq -\frac{\alpha_\gn^{f,g,2}}\gamma\ell(x,u) + \frac{1-\gamma}\gamma V_\gn^f(x))
    \end{align*}
    which shows \eqref{eq:rDPg}.
    Applying the same procedure starting with \eqref{eq:rDPn} yields \eqref{eq:rDPgn}.

\noindent\textbf{Proof of item 1) of Theorem \ref{thm:main}.}
    Let $\overline c\in\R_{\geq0}\cup\{\infty\}$ such that $\mc L_\gn^f(\mc S) = \{x\in\mc S~|~V_\gn^f(x)\leq\overline c\}$.
    Let $x_k,k\in\N_{0}$ be a solution to \eqref{eq:closed_loop} with $x_0\in\mc L_\gn^f(\mc S)$.
    We will prove by induction that $x_k\in\mc L_\gn^f(\mc S)$ for all $k\in\N_0$.
    %If $\mc S=\R^n$, this is trivial since then $\mc L_\gn^f(\mc S)=\mc S$.
    %Now assume $\aggnn\geq0$, then we conduct a proof by induction.
    The base case is already established by choice of $x_0$.
    Now suppose that $x_k\in\mc L_\gn^f(\mc S)$ for some $k\in\N_{0}$ and let $u_k\in\mc U_\gn^f(x_k)$ such that $x_{k+1} = g(x_k,u_k)$.
    Since $x_k\in\mc L_\gn^f(\mc S)\subseteq\mc S$, we can apply inequality \eqref{eq:rDPgn} for $x_k$ and $u_k$, by which 
    \begin{align}
        &V_\gn^f(x_{k+1}) = V_\gn^f(g(x_k,u_k))\nonumber\\
        &\leq V_\gn^f(x_k) - \aggnn\ell(x_k,u_k)\nonumber\\
        &\leq V_\gn^f(x_k) - \aggnn\|x_k\|_Q^2\leq \left(1- \tfrac\aggnn B\right)V_\gn^f(x_k)\label{eq:proof:exponential:1}\\
        &\leq V_\gn^f(x_k) \leq \overline c\label{eq:proof:exponential:2}
    \end{align}
    where we used \eqref{eq:B-cost-controllable} and $\aggnn\geq0$.
    Thus, $x_{k+1}\in\mc L_\gn^f(\mc S)$, completing the proof by induction. Iterating \eqref{eq:proof:exponential:1},
        $\|x_k\|_Q^2\leq V_\gn^f(x_k) \leq \left(1\!-\!\frac\aggnn B\right)^kV_\gn^f(x_0) \leq B\left(1\!-\!\frac\aggnn B\right)^k\|x_0\|_Q^2$
    for every $k\in\N_{0}$.
    Taking square roots completes the proof.

\noindent\textbf{Proof of item 2) of Theorem \ref{thm:main}.}
    Let $x_k$ and $u_k$, $k\in\N_0$ be as in the proof of item 1), by which $x_k\in\mc S$ for all $k\in\N_0$ and we can apply \eqref{eq:rDPgn}.
    Rearranging and summing yields a telescoping sum, whereby $\aggnn\sum_{k=0}^\infty\ell(x_k,u_k) \leq \sum_{k=0}^\infty \left(V_\gn^f(x_k) - V_\gn^f(x_{k+1}) \right)\leq V_\gn^f(x_0)$.
    Since this holds for all solutions of \eqref{eq:closed_loop}, inequality \eqref{eq:thm:suboptimality2} follows by the definition of $\overline J_\oi^g(x,\mc U_\gn^f)$ in \eqref{eq:performance2}.
    Using \eqref{eq:rDPn} instead and multiplying by $\gamma^k$ yields
    $\agnn\sum_{k=0}^\infty\gamma^k\ell(x_k,u_k) \leq \sum_{k=0}^\infty\gamma^k \left(V_\gn^f(x_k) - \gamma V_\gn^f(x_{k+1}) \right)\leq V_\gn^f(x_0),$
   which proves \eqref{eq:thm:suboptimality} in similar fashion.
   
\noindent\textbf{Proof of item 3) of Theorem \ref{thm:main}.}
Consider a sequence of states $x_k$ and controls $u_k$, where $k\in\{0,\dots,N-1\}$, constructed by $x_0:=x$ and $u_k\in\mc U_{\gamma,N-k}^f(x_k)$, $x_{k+1}:=g(x_k,u_k)$ for $k\in\{0,\dots,N-1\}$.
We will show that $x_k\in\mc S$ for all $k\in\{0,\dots,N-1\}$.
In case (a) this is trivial as $\mc S=\R^n$.

\underline{If (b) holds:} We will show by induction that $V_{\gamma,N-k}^f(x_k)\leq\overline c$ for each $k\in\{0,\dots,N-1\}$.
For $k=0$, $V_{\gamma,N}^f(x_0)\leq\overline c$ holds by choice of $x=x_0$.
Now suppose $V_{\gamma,N-k}^f(x_k)\leq\overline c$ for some $k\in\{0,\dots,N-2\}$, then $\|x_k\|_Q^2\leq V_{\gamma,N-k}^f(x_k)\leq\overline c$, and hence $x_k\in\mc S$ by definition of $\overline c$.
Thus, we can apply \eqref{eq:rDPg} to $x_k$ and $u_k$, which yields $V_{\gamma,N-k-1}^f(x_{k+1})-V_{\gamma,N-k}^f(x_k) \leq -\alpha_{\underline{\smash{\gamma}},N-k}\ell(x_k,u_k) \leq -\aggn\ell(x_k,u_k)\leq0$.
In particular, $V_{\gamma,N-k-1}^f(x_{k+1}) \leq V_{\gamma,N-k}^f(x_k) \leq \overline c$ and $x_{k+1}\in\mc S$, completing the proof by induction.

\underline{If (c)  holds:} We begin with the following two lemmas.
\begin{lem}\label{lem:different_horizons}
    For any $m\in\{1,\dots,N-1\}$ and $x\in\R^n$,
    \begin{align}
        V_\gn^f(x)\leq \left(1+B(1-\tfrac1B)^m\right)V_{\gamma,m}^f(x).
    \end{align}
\end{lem}
\begin{proof}
    Let $\mathbf u_{m}\in\mathbb U^m$ such that $V_{\gamma,m}^f(x) = J_{\gamma,m}^f(x,\mathbf u_m)$, and let $x_k:=\varphi^f_x(k,\mathbf u_{k})$ for $k\in\{0,\dots,m-1\}$.
    By Proposition~\ref{prop:1.bellman} and SA\ref{assm:cost_controllable}, $\gamma V_{\gamma,m-k-1}^f(x_{k+1}) = V_{\gamma,m-k}^f(x_{k})-\ell(x_k,u_k)\leq(1-\tfrac1B)V_{\gamma,m-k}^f(x_k)$ for all $k\in\{0,\dots,m-2\}$.
    Iterating this yields $\|x_{m-1}\|_Q^2 = V_{\gamma,1}^f(x_{m-1})\leq(\tfrac1\gamma(1-\frac1B))^{m-1}V_{\gamma,m}^f(x)$.
    By Proposition \ref{prop:1.bellman} and then using SA\ref{assm:cost_controllable} and this previous inequality,
    \begin{align}
        V_\gn^f(x) &\leq J_{\gamma,m-1}^f(x,\mathbf u_{m-1})+\gamma^{m-1}V_{\gamma,N-m+1}^f(x_{m-1})\nonumber\\
        &\leq V_{\gamma,m}^f(x)-\gamma^{m-1}\|x_{m-1}\|_Q^2 + \gamma^{m-1}B\|x_{m-1}\|_Q^2\nonumber\\
        &\leq V_{\gamma,m}^f(x) + (B-1)(1-\tfrac1B)^{m-1}V_{\gamma,m}^f(x)\nonumber\\
        &=(1+B(1-\tfrac1B)^m)V_{\gamma,m}^f(x),
    \end{align}
    which completes the proof of the lemma.
\end{proof}

\begin{lem}\label{lem:technical_inequality}
    For any $B\geq3$, $N\geq2B\log(B)$ and $k\in\{1,\dots,N-2\}$, it holds that
    %\begin{align*}
        $\mu(k):=\left(1+B\left(1-\tfrac1B\right)^{N-k-1}\right)\left(1-\tfrac{2\log(B)}{N}\right)^{k+1}\leq1.$
    %\end{align*}
\end{lem}
%The proof of Lemma \ref{lem:technical_inequality} can be found in \cite{moldenhauer2026discounted}.
\begin{proof}
    Let $a:=1-\tfrac1B$ and $c:=1-\tfrac{2\log(B)}{N}$.
    Then $\mu(k) = c^{k+1} + Ba^{N}(\tfrac ca)^{k+1}$, which is the sum of two terms that are convex in $k$ (since $c,\tfrac ca\geq0$) and thereby convex in $k$.
    Hence, $\mu(k)\leq\max\{\mu(1),\mu(N-2)\}$ for all $k\in\{1,\dots,N-2\}$,
    which means that we only need to verify $\mu(k)\leq1$ for $k=1$ and $k=N-2$.\\
    \underline{$k=1$:} Using the basic inequality $1+t\leq e^t$ for $t=Ba^{N-2}$ and $t=-2\log(B)/N$, we get $\mu(1)=(1+Ba^{N-2})(1-{2\log(B)}/{N})^2\leq\exp(Ba^{N-2}-{4\log(B)}/N)$.
    We will show $Ba^{N-2}-{4\log(B)}/N\leq0$, which would imply $\mu(1)\leq1$.
    Using $a=1-\frac1B\leq e^{-1/B}$ and taking logarithms, we get that it is sufficient to show $\log(B)-(N-2)/B\leq \log(4\log(B))-\log(N)$.
    Rearranged, this becomes $\nu(N):=\log(4\log(B))-\log(B)+(N-2)/B-\log(N)\geq0$.
    %Since $1\geq\aggn\geq 2B\log(B)/N$ by condition (c), it must necessarily hold that $N\geq 2B\log(B)$.
    Treating $\nu$ as a function on $\R_{\geq0}$, $\tfrac{d}{dt}\nu(t)=1/B-1/t\geq0$ for all $t\geq 2B\log(B)$, since $2B\log(B)\geq B$ for $B\geq3$.
    With that and $N\geq2B\log(B)$, $\nu(N)\geq \nu(2B\log(B))=\log(2)-2/B$, which is positive for $B\geq3$.
    This completes the proof for $k=1$.\\
    \underline{$k=N-2$.} As in the proof for $k=1$ and since $B\geq3$, we have $N\geq 2B\log(B)\geq 2$, and hence $2\log(B)/N\geq \log(B)/(N-1)$.
    By this and $1-t\leq e^{-t}$ for $t=\log(B)/(N-1)$, we get
    $\mu(N-2)=B(1-2\log(B)/N)^{N-1}\leq B(1-\log(B)/(N-1))^{N-1}\leq Be^{-\log(B)}=1$.
\end{proof}
We now proceed with the proof assuming (c). Let $\overline c\in\R_{\geq0}\cup\{\infty\}$ such that $\mc L_\gn^f(\mc S) = \{x\in\R^n~|~V_\gn^f(x)\leq\overline c\}$.
We will show by induction that $V_{\gamma,N}^f(x_k)\leq\overline c$ for each $k\in\{0,\dots,N-1\}$.
For $k=0$, $V_{\gamma,N}^f(x_0)\leq\overline c$ holds by choice of $x=x_0$.
For $k=1$, since $x_0\in\mc S$ and $u_0\in\mc U_\gn^f(x_0)$, we can apply  \eqref{eq:rDPgn}, which yields $V_\gn^f(x_1)\leq V_\gn^f(x_0)-\aggnn\ell(x_0,u_0)\leq V_\gn^f(x_0)$ since $\aggnn\geq0$.
%(this can be shown using $\aggn\geq2B\log(B+1)/N, \gamma\geq1-1/B$ and $B\geq2$; we skip the details),
%we have $V_\gn^f(x_1)\leq\overline c$.
Now suppose $V_{\gamma,N}^f(x_k)\leq\overline c$ for some $k\in\{1,\dots,N-2\}$, then $x_k\in\mc S$ by definition of $\overline c$.
Thus, we can apply \eqref{eq:rDPg} to $x_k$ and $u_k$, which yields $V_{\gamma,N-k-1}^f(x_{k+1})\leq V_{\gamma,N-k}^f(x_k) -\alpha_{\underline{\smash{\gamma}},N-k}\ell(x_k,u_k) \leq V_{\gamma,N-k}^f(x_k)-\aggn\ell(x_k,u_k)\leq(1-\frac\aggn B)V_{\gamma,N-k}^f(x_k)$.
Iterating this over all previous $k$ yields $V_{\gamma,N-k-1}^f(x_{k+1})\leq (1-\frac\aggn B)^{k+1}V_\gn^f(x_0)$.
Combining this with Lemma \ref{lem:different_horizons}, $\aggn\geq2B\log(B)/N$ and Lemma~\ref{lem:technical_inequality} (which is applicable since $N\geq2B\log(B)$ by $1\geq\aggn\geq2B\log(B)/N$), we have $V_{\gamma,N}^f(x_{k+1})\leq (1+B(1-\frac1B)^{N-k-1})(1-\frac\aggn B)^{k+1}V_\gn^f(x_0) \leq V_\gn^f(x_0)\leq\overline c$ and $x_{k+1}\in\mc S$, completing the induction step for (c).
%We now know that $x_k\in\mc S$ for all $k\in\{0,\dots,N-1\}$.

In any case, because $x_k\in\mc S$ and $u_k\in\mc U_{\gamma,N-k}^f(x_k)$ for all $k\in\{0,\dots,N-1\}$, we can apply \eqref{eq:rDP0} to obtain $\gamma V_{\gamma,N-k-1}^f(x_{k+1})-V_{\gamma,N-k}^f(x_k)\leq-\alpha_{\gamma,N-k}\ell(x_k,u_k) \leq -\alpha_\gn\ell(x_k,u_k)$ for $k\leq N-2$, which also trivially holds for $k=N-1$. % since $V_{\gamma,0}^f:\equiv0$, $V_{\gamma,1}^f(x_{N-1})=\ell(x_{N-1},u_{N-1})$ and $\alpha_{\gamma,1}\leq1$.
%In any case, because $x_k\in\mc S$ and $u_k\in\mc U_{\gamma,N-k}^f(x_k)$ for all $k\in\{0,\dots,N-1\}$, we can apply \eqref{eq:rDP0} to obtain $\gamma V_{\gamma,N-k-1}(x_{k+1})-V_{\gamma,N-k}(x_k)\leq-\alpha_{\gamma,N-k}\ell(x_k,u_k) \leq -\alpha_\gn\ell(x_k,u_k)$.
Adding these up multiplied with $\gamma^k$ 
%and using the definition of $V_\gn^g(x)$
yields $\agn V_\gn^g(x) \leq\sum_{k=0}^{N-1}\gamma^k\agn\ell(x_k,u_k)
    \leq \sum_{k=0}^{N-1}\gamma^k \left(V_{\gamma,N-k}^f(x_k) - \gamma V_{\gamma,N-k-1}^f(x_{k+1}) \right)\leq V_\gn^f(x),$
% \begin{align}
%     \agn V_\gn^g(x) &\leq\sum_{k=0}^{N-1}\gamma^k\agn\ell(x_k,u_k)\nonumber\\
%     &\leq \sum_{k=0}^{N-1}\gamma^k \left(V_{\gamma,N-k}^f(x_k) - \gamma V_{\gamma,N-k-1}^f(x_{k+1}) \right)\nonumber\\
%     &\leq V_\gn^f(x),
% \end{align}
which completes the proof of item 3).

\end{document}